\documentclass[12pt]{article}%
\usepackage{amsfonts}
\usepackage{sw20bams}
\usepackage{amsmath}
\usepackage{amssymb}
\usepackage{graphicx}%
\providecommand{\U}[1]{\protect\rule{.1in}{.1in}}
\begin{document}

\title{M/G/1-FIFO Queue with Uniform Service Times}
\author{Steven Finch}
\date{October 18, 2022}
\maketitle

\begin{abstract}
An exact formula for the equilibrium M/U/1 waiting time density is now
effectively known. What began as a numeric exploration became a symbolic
banquet. \ Inverse Laplace transforms provided breadcrumbs in the trail; delay
differential equations subsequently gave clear-cut precision. \ We also remark
on tail probability asymptotics and on queue lengths.

\end{abstract}

\footnotetext{Copyright \copyright \ 2022 by Steven R. Finch. All rights
reserved.}Consider a first-in-first-out M/G/1 queue alongside unlimited
waiting space, where the input process is Poisson with rate $\lambda$ and the
service times are independent Uniform[$a,b$] random variables with mean
$1/\mu=(a+b)/2$. \ Let $W_{\text{que}}$ denote the waiting time in the queue
(prior to service). \ Under equilibrium (steady-state) conditions and traffic
intensity (load) $\rho=\lambda/\mu<1$, the probability density function $f(x)$
of $W_{\text{que}}$ has Laplace transform \cite{STGH-que, Neu-que}%
\[
F(s)=\lim_{\varepsilon\rightarrow0^{+}}%
{\displaystyle\int\limits_{-\varepsilon}^{\infty}}
\exp(-s\,x)f(x)dx=\frac{(1-\rho)s}{s-\lambda+\lambda\,\Theta(s)}%
=F_{\text{alt}}(s)+1-\rho
\]
where%
\[%
\begin{array}
[c]{ccc}%
\theta(x)=\left\{
\begin{array}
[c]{lll}%
1/(b-a) &  & a\leq x\leq b\\
0 &  & \text{otherwise}%
\end{array}
\right.  , &  & \Theta(s)=\dfrac{\exp(-a\,s)-\exp(-b\,s)}{(b-a)s}.
\end{array}
\]
From%
\[
(1-\rho)s=s\,F(s)-\lambda\,F(s)\left[  1-\Theta(s)\right]
\]
we have%
\[%
\begin{array}
[c]{ccccc}%
(1-\rho)s+\lambda\,F(s)\left[  1-\Theta(s)\right]  =s\,F(s), &  & \text{i.e.,}
&  & F(s)=1-\rho+\lambda\,F(s)\left[  \dfrac{1}{s}-\dfrac{\Theta(s)}%
{s}\right]
\end{array}
\]
hence%
\[
f(x)=(1-\rho)\delta(x)+\kappa+\lambda%
{\displaystyle\int\limits_{0}^{x}}
f(t)\left[  1-%
{\displaystyle\int\limits_{0}^{x-t}}
\theta(u)du\right]  dt
\]
where $\delta(x)$ is the Dirac delta and $\kappa=\rho(\mu-\lambda)$.
\ Differentiating, we obtain
\begin{align*}
f^{\prime}(x) &  =\lambda f(x)\left[  1-0\right]  +\lambda%
{\displaystyle\int\limits_{0}^{x}}
f(t)\left[  0-\theta(x-t)\right]  dt\\
&  =\lambda f(x)-\lambda%
{\displaystyle\int\limits_{\max\{x-b,0\}}^{x-a}}
f(t)\frac{1}{b-a}dt.
\end{align*}
There are three cases, to be examined separately. \ For simplicity, we set
$\lambda=2$, $\mu=3$ and choose $a$, $b$ appropriately. \ However, it can be
shown (in general) that%
\[%
\begin{array}
[c]{ccc}%
f^{\prime}(x)=\lambda f(x), &  & \lim\limits_{\varepsilon\rightarrow0^{+}%
}f(\varepsilon)=\kappa
\end{array}
\]
is valid for $0<x<a\leq b$, i.e., $f(x)$ is equal to $\kappa\exp(\lambda\,x)$;
and%
\[%
\begin{array}
[c]{ccccc}%
f^{\prime\prime}(x)=\lambda f^{\prime}(x)-\dfrac{\lambda\,\mu}{2}f(x), &  &
f(0^{+})=\kappa, &  & f^{\prime}(0^{+})=\kappa\left(  \lambda-\dfrac{\mu}%
{2}\right)
\end{array}
\]
is valid for $0=a<x<b$ (since $b/2=1/\mu$), i.e., $f(x)$ is equal to
\[
\kappa\exp\left(  \frac{\lambda\,x}{2}\right)  \left[  \cos\left(  \frac{1}%
{2}\sqrt{\lambda(2\mu-\lambda)}x\right)  -\frac{\mu-\lambda}{\sqrt
{\lambda(2\mu-\lambda)}}\sin\left(  \frac{1}{2}\sqrt{\lambda(2\mu-\lambda
)}x\right)  \right]  .
\]
The indicated conditions are true due to formulas%
\[%
\begin{array}
[c]{ccc}%
f(0^{+})=\lim\limits_{s\rightarrow1\cdot\infty}s\,F_{\text{alt}}(s), &  &
f^{\prime}(0^{+})=\lim\limits_{s\rightarrow1\cdot\infty}\left[  s^{2}%
F_{\text{alt}}(s)-s\,f(0^{-})\right]
\end{array}
\]
that proceed from the initial value theorem \cite{LMT-que, Sch-que, Hen-que}.

\section{Case One ($0<a<b$)}

Set $a=\frac{1}{12}$, $b=\frac{7}{12}$ and%
\[%
\begin{array}
[c]{ccc}%
c_{n}=\left\{
\begin{array}
[c]{ccc}%
-\infty &  & \text{if }n=0\\
\frac{n}{12} &  & \text{if }n\geq1
\end{array}
\right.  , &  & f_{0}(x)=\frac{2}{3}\exp(2\,x),
\end{array}
\]%
\[%
\begin{array}
[c]{ccccc}%
i_{n}(x)=4%
{\displaystyle\int\limits_{c_{n}}^{x-\frac{1}{12}}}
f_{n}(t)dt, &  & j_{n}=4%
{\displaystyle\int\limits_{c_{n}}^{c_{n+1}}}
f_{n}(t)dt, &  & k_{n}(x)=4%
{\displaystyle\int\limits_{x-\frac{7}{12}}^{c_{n+1}}}
f_{n}(t)dt.
\end{array}
\]
A sequence of functions is defined iteratively as follows:%
\[%
\begin{array}
[c]{ccc}%
f_{n}^{\prime}(x)=2f_{n}(x)-i_{n-1}(x)-%
{\displaystyle\sum\limits_{p=\max\{0,n-6\}}^{n-2}}
j_{p}+k_{n-7}(x), &  & f_{n}\left(  \tfrac{n}{12}\right)  =f_{n-1}\left(
\tfrac{n}{12}\right)
\end{array}
\]
where the empty sum convention holds for $n=1$ and it is understood that
$k_{q}=0$ for $q<0$. \ Since $f_{n}(x)$ is a degree $n$ polynomial in $x$ with
coefficients of the form $R\exp(r+2x)$, where $R$, $r$ are rational numbers,
the integrals $i_{n}$, $j_{n}$, $k_{n}$ all possess closed-form expressions.
\ Therefore the differential equation for $f_{n+1}(x)$ can be solved exactly.
\ Stitching the fragments together gives rise to the density function
\[
f(x)=f_{\left\lfloor 12x\right\rfloor }(x)
\]
pictured in Figure 1, where $\left\lfloor y\right\rfloor $ denotes the
greatest integer $\leq y$.

Let us illustrate in greater detail:\
\[
i_{0}(x)=4%
{\displaystyle\int\limits_{-\infty}^{x-\frac{1}{12}}}
\tfrac{2}{3}\exp(2\,t)dt=\tfrac{4}{3}\exp\left(  -\tfrac{1}{6}+2x\right)
\]
thus%
\[%
\begin{array}
[c]{ccc}%
f_{1}^{\prime}(x)=2f_{1}(x)-i_{0}(x), &  & f_{1}\left(  \tfrac{1}{12}\right)
=f_{0}\left(  \tfrac{1}{12}\right)  =\tfrac{2}{3}\exp(\tfrac{1}{6})
\end{array}
\]
implying%
\[
f_{1}(x)=\tfrac{2}{3}\exp\left(  2x\right)  +\tfrac{1}{9}\exp\left(
-\tfrac{1}{6}+2x\right)  -\tfrac{4}{3}\exp\left(  -\tfrac{1}{6}+2x\right)  x.
\]
Continuing:
\begin{align*}
i_{1}(x)  &  =4%
{\displaystyle\int\limits_{\frac{1}{12}}^{x-\frac{1}{12}}}
\left[  \tfrac{2}{3}\exp\left(  2t\right)  +\tfrac{1}{9}\exp\left(  -\tfrac
{1}{6}+2t\right)  -\tfrac{4}{3}\exp\left(  -\tfrac{1}{6}+2t\right)  t\right]
dt\\
&  =-\tfrac{4}{3}-\tfrac{4}{3}\exp(\tfrac{1}{6})+\tfrac{16}{9}\exp\left(
-\tfrac{1}{3}+2x\right)  +\tfrac{4}{3}\exp\left(  -\tfrac{1}{6}+2x\right)
-\tfrac{8}{3}\exp\left(  -\tfrac{1}{3}+2x\right)  x,
\end{align*}%
\[
j_{0}=4%
{\displaystyle\int\limits_{-\infty}^{\frac{1}{12}}}
\tfrac{2}{3}\exp(2\,t)dt=\tfrac{4}{3}\exp(\tfrac{1}{6})
\]
thus%
\[%
\begin{array}
[c]{ccc}%
f_{2}^{\prime}(x)=2f_{2}(x)-i_{1}(x)-j_{0}, &  & f_{2}\left(  \tfrac{1}%
{6}\right)  =f_{1}\left(  \tfrac{1}{6}\right)  =-\tfrac{1}{9}\exp(\tfrac{1}%
{6})+\tfrac{2}{3}\exp(\tfrac{1}{3})
\end{array}
\]
implying%
\begin{align*}
f_{2}(x)  &  =-\tfrac{2}{3}+\tfrac{2}{3}\exp(2x)+\tfrac{25}{27}\exp\left(
-\tfrac{1}{3}+2x\right)  +\tfrac{1}{9}\exp\left(  -\tfrac{1}{6}+2x\right) \\
&  -\tfrac{16}{9}\exp\left(  -\tfrac{1}{3}+2x\right)  x-\tfrac{4}{3}%
\exp\left(  -\tfrac{1}{6}+2x\right)  x+\tfrac{4}{3}\exp\left(  -\tfrac{1}%
{3}+2x\right)  x^{2}.
\end{align*}
The pattern is maintained:
\[
f_{3}^{\prime}(x)=2f_{3}(x)-i_{2}(x)-j_{1}-j_{0},
\]%
\[
f_{4}^{\prime}(x)=2f_{4}(x)-i_{3}(x)-j_{2}-j_{1}-j_{0},
\]%
\[
f_{5}^{\prime}(x)=2f_{5}(x)-i_{4}(x)-j_{3}-j_{2}-j_{1}-j_{0},
\]%
\[
f_{6}^{\prime}(x)=2f_{6}(x)-i_{5}(x)-j_{4}-j_{3}-j_{2}-j_{1}-j_{0}%
\]
until $x>\frac{7}{12}$ (the equation length becomes fixed and $k$ replaces the
rightmost $j$):%
\[
f_{7}^{\prime}(x)=2f_{7}(x)-i_{6}(x)-j_{5}-j_{4}-j_{3}-j_{2}-j_{1}-k_{0}(x),
\]%
\[
f_{8}^{\prime}(x)=2f_{8}(x)-i_{7}(x)-j_{6}-j_{5}-j_{4}-j_{3}-j_{2}-k_{1}(x),
\]%
\[
f_{9}^{\prime}(x)=2f_{9}(x)-i_{8}(x)-j_{7}-j_{6}-j_{5}-j_{4}-j_{3}-k_{2}(x),
\]%
\[
f_{10}^{\prime}(x)=2f_{10}(x)-i_{9}(x)-j_{8}-j_{7}-j_{6}-j_{5}-j_{4}%
-k_{3}(x).
\]
We count $n(n+5)/2$ terms in each $f_{n}(x)$ upon expansion, at least for
$1\leq n\leq35$. \ 

When numerically evaluating large symbolic expressions, it is important to
ensure that the working precision of floating point quantities is suitably
high. To evaluate $f(2)$ may require several hundred decimal digits because,
for instance, the first two of the 348 numerator terms comprising $f_{24}(x)$
are%
\begin{align*}
&  -31343712612206064875238458599056650210472221756256360\,\exp(4)\\
&  +1235688973308606091819588575256480179309880006812893184\,\exp(2x)\\
&  \approx(-1.7113\times10^{54})+(6.7466\times10^{55})
\end{align*}
upon setting $x=2$. \ The subtraction of nearly equal numbers, such as these,
will lead to a horrific loss of precision unless appropriate care is taken.

\section{Case Two ($0=a<b$)}

Set $a=0$, $b=\frac{2}{3}$ and%
\[
g_{0}(x)=\tfrac{1}{6}\exp(x)\left[  4\cos\left(  \sqrt{2}x\right)  -\sqrt
{2}\sin\left(  \sqrt{2}x\right)  \right]  .
\]
A sequence of functions is defined iteratively as follows:%
\[
g_{n}^{\prime\prime}(x)=2g_{n}^{\prime}(x)-3g_{n}(x)+3g_{n-1}\left(
x-\tfrac{2}{3}\right)  ,
\]%
\[%
\begin{array}
[c]{ccc}%
g_{n}\left(  \tfrac{2n}{3}\right)  =g_{n-1}\left(  \tfrac{2n}{3}\right)  , &
& g_{n}^{\prime}\left(  \tfrac{2n}{3}\right)  =\delta_{n-1}+g_{n-1}^{\prime
}\left(  \tfrac{2n}{3}\right)
\end{array}
\]
where $\delta_{m}$ is the Kronecker delta. \ Stitching the fragments together
gives rise to
\[
g(x)=g_{\left\lfloor \frac{3x}{2}\right\rfloor }(x)
\]
pictured in Figure 2.

Let us illustrate in greater detail:
\[
g_{0}^{\prime}(x)=\tfrac{1}{6}\exp(x)\left[  2\cos\left(  \sqrt{2}x\right)
-5\sqrt{2}\sin\left(  \sqrt{2}x\right)  \right]
\]
thus%
\[
g_{1}^{\prime\prime}(x)=2g_{1}^{\prime}(x)-3g_{1}(x)+3g_{0}\left(  x-\tfrac
{2}{3}\right)  ,
\]%
\[%
\begin{array}
[c]{ccc}%
g_{1}\left(  \tfrac{2}{3}\right)  =g_{0}\left(  \tfrac{2}{3}\right)  , &  &
g_{1}^{\prime}\left(  \tfrac{2}{3}\right)  =1+g_{0}^{\prime}\left(  \tfrac
{2}{3}\right)
\end{array}
\]
implying
\begin{align*}
g_{1}(x)  &  =\tfrac{1}{6}\exp(x)\left[  4\cos\left(  \sqrt{2}x\right)
-\sqrt{2}\sin\left(  \sqrt{2}x\right)  \right] \\
&  -\tfrac{1}{24}\exp\left(  -\tfrac{2}{3}+x\right)  \left[  4\cos\left(
\sqrt{2}\left(  -\tfrac{2}{3}+x\right)  \right)  -\sqrt{2}\sin\left(  \sqrt
{2}\left(  -\tfrac{2}{3}+x\right)  \right)  \right] \\
&  +\tfrac{1}{4}\exp\left(  -\tfrac{2}{3}+x\right)  \left[  \cos\left(
\sqrt{2}\left(  -\tfrac{2}{3}+x\right)  \right)  +2\sqrt{2}\sin\left(
\sqrt{2}\left(  -\tfrac{2}{3}+x\right)  \right)  \right]  x.
\end{align*}
Beyond this point, the derivatives at $x=\frac{2n}{3}$ match: \
\begin{align*}
g_{1}^{\prime}(x)  &  =\tfrac{1}{6}\exp(x)\left[  2\cos\left(  \sqrt
{2}x\right)  -5\sqrt{2}\sin\left(  \sqrt{2}x\right)  \right] \\
&  +\tfrac{1}{24}\exp\left(  -\tfrac{2}{3}+x\right)  \left[  4\cos\left(
\sqrt{2}\left(  -\tfrac{2}{3}+x\right)  \right)  +17\sqrt{2}\sin\left(
\sqrt{2}\left(  -\tfrac{2}{3}+x\right)  \right)  \right] \\
&  +\tfrac{1}{4}\exp\left(  -\tfrac{2}{3}+x\right)  \left[  5\cos\left(
\sqrt{2}\left(  -\tfrac{2}{3}+x\right)  \right)  +\sqrt{2}\sin\left(  \sqrt
{2}\left(  -\tfrac{2}{3}+x\right)  \right)  \right]  x
\end{align*}
thus%
\[
g_{2}^{\prime\prime}(x)=2g_{2}^{\prime}(x)-3g_{2}(x)+3g_{1}\left(  x-\tfrac
{2}{3}\right)  ,
\]%
\[%
\begin{array}
[c]{ccc}%
g_{2}\left(  \tfrac{4}{3}\right)  =g_{1}\left(  \tfrac{4}{3}\right)  , &  &
g_{2}^{\prime}\left(  \tfrac{4}{3}\right)  =g_{1}^{\prime}\left(  \tfrac{4}%
{3}\right)
\end{array}
\]
implying%
\begin{align*}
g_{2}(x)  &  =\tfrac{1}{6}\exp(x)\left[  4\cos\left(  \sqrt{2}x\right)
-\sqrt{2}\sin\left(  \sqrt{2}x\right)  \right] \\
&  -\tfrac{1}{24}\exp\left(  -\tfrac{2}{3}+x\right)  \left[  4\cos\left(
\sqrt{2}\left(  -\tfrac{2}{3}+x\right)  \right)  -\sqrt{2}\sin\left(  \sqrt
{2}\left(  -\tfrac{2}{3}+x\right)  \right)  \right] \\
&  +\tfrac{1}{4}\exp\left(  -\tfrac{2}{3}+x\right)  \left[  \cos\left(
\sqrt{2}\left(  -\tfrac{2}{3}+x\right)  \right)  +2\sqrt{2}\sin\left(
\sqrt{2}\left(  -\tfrac{2}{3}+x\right)  \right)  \right]  x\\
&  -\tfrac{1}{192}\exp\left(  -\tfrac{4}{3}+x\right)  \left[  8\cos\left(
\sqrt{2}\left(  -\tfrac{4}{3}+x\right)  \right)  -29\sqrt{2}\sin\left(
\sqrt{2}\left(  -\tfrac{4}{3}+x\right)  \right)  \right] \\
&  +\tfrac{1}{32}\exp\left(  -\tfrac{4}{3}+x\right)  \left[  17\cos\left(
\sqrt{2}\left(  -\tfrac{4}{3}+x\right)  \right)  -2\sqrt{2}\sin\left(
\sqrt{2}\left(  -\tfrac{4}{3}+x\right)  \right)  \right]  x\\
&  -\tfrac{3}{32}\exp\left(  -\tfrac{4}{3}+x\right)  \left[  4\cos\left(
\sqrt{2}\left(  -\tfrac{4}{3}+x\right)  \right)  -\sqrt{2}\sin\left(  \sqrt
{2}\left(  -\tfrac{4}{3}+x\right)  \right)  \right]  x^{2}.
\end{align*}

\section{Case Three ($0<a=b$)}

Set $a=b=\frac{1}{3}$ and%
\[
h_{0}(x)=\tfrac{2}{3}\exp(2\,x).
\]
A sequence of functions is defined iteratively as follows:%
\[
h_{n}^{\prime}(x)=2h_{n}(x)-2h_{n-1}\left(  x-\tfrac{1}{3}\right)  ,
\]%
\[
h_{n}\left(  \tfrac{n}{3}\right)  =-\tfrac{2}{3}\delta_{n-1}+h_{n-1}\left(
\tfrac{n}{3}\right)  .
\]
The role of $\delta_{m}$ here is more pronounced than in Section 2: a jump
discontinuity occurs in the density at $x=\frac{1}{3}$ (as opposed to merely a
sharp corner) . \ Stitching the fragments together gives rise to%
\[
h(x)=h_{\left\lfloor 3x\right\rfloor }(x)
\]
pictured in Figure 3. \ 

We could do as before, indicating steps leading to $h_{1}(x)$ and $h_{2}(x)$.
A\ classical result due to Erlang \cite{Erl-que, IvS-que, Pra-que, Tij1-que}:%
\[
h_{n}(x)=\left(  1-\rho\right)  \frac{d}{dx}\left\{
{\displaystyle\sum\limits_{m=0}^{n}}
(-1)^{m}\frac{\lambda^{m}\left(  x-a\,m\right)  ^{m}}{m!}\exp\left[
\lambda(x-a\,m)\right]  \right\}
\]
renders this listing unnecessary (with $1-\rho=\frac{1}{3}$, $\lambda=2$,
$a=\frac{1}{3}$). \ We wonder if such a formula (for the M/D/1 queue)
possesses an analog for the M/U/1 queue.

\section{Personal Notes}

I\ taught a semester-long class in statistical programming at Harvard
University (as a preceptor) for nearly ten years. \ A favorite set of problems
began with an M/M/1 example \cite{HL-que} involving a hospital emergency room (ER).

\begin{quotation}
\noindent Patients (clients) arrive according to a Poisson process with rate
$\lambda$. \ One\ doctor (server)\ is available to treat them. \ The doctor,
when busy, treats patients with rate $\mu$. \ More precisely: interarrival
times are exponentially distributed with mean $1/\lambda$ and treatment
lengths are exponentially distributed with mean $1/\mu$. \ The ER is open 24
hours a day, 7 days a week. \ Patients must wait until the doctor is free and
are treated in the order via which they arrive. \ Simulate the performance of
the ER over many weeks. \ What can be said about waiting time in queue per
patient (excluding service time)? \ Determine the mean, variance, mode and
median of $W_{\text{que}}$ to as high accuracy as possible. \ Assume for this
purpose that the expected number of arriving patients per hour is $\lambda=2$
and that the expected number of treatment completions per hour (for a
continuously busy ER) is $\mu=3$.
\end{quotation}

\noindent Putting aside experiment in favor of pristine theory, the Laplace
transform of service density $\mu\exp(-\mu\,x)$ is $\mu/(\mu+s)$. \ Formulas
for the mean and variance of $W_{\text{que}}$ follow immediately from%
\[
F(s)=\frac{(1-\rho)s}{s-\lambda+\lambda\frac{\mu}{\mu+s}};
\]
consequently%
\[%
\begin{array}
[c]{ccc}%
\text{mean}=-F^{\prime}(0)=\dfrac{\lambda}{\mu(\mu-\lambda)}, &  &
\text{variance}=F^{\prime\prime}(0)-F^{\prime}(0)^{2}=\dfrac{\lambda
(2\mu-\lambda)}{\mu^{2}(\mu-\lambda)^{2}}%
\end{array}
\]
giving $\frac{2}{3}$ and $\frac{8}{9}$ respectively. \ An expression for the
density of $W_{\text{que}}$:
\[
f(x)=(1-\rho)\delta(x)+\rho(\mu-\lambda)\exp(-(\mu-\lambda)x)
\]
shows that the mode is $0$; integrating and solving the equation%
\[
1-\rho-\rho\left[  \exp(-(\mu-\lambda)x)-1\right]  =\dfrac{1}{2}%
\]
implies%
\[%
\begin{array}
[c]{ccccc}%
\dfrac{1}{2}=\rho\exp(-(\mu-\lambda)x), &  & \text{i.e.,} &  & x=\dfrac{1}%
{\mu-\lambda}\ln\left(  \dfrac{2\lambda}{\mu}\right)
\end{array}
\]
giving the median to be $\ln(4/3)=0.28768207...$.

The aforementioned problem set continued with an M/G/1 example, involving the
same ER parameters, but with Uniform[$a,b$] treatment lengths. \ Since
$1/3=(a+b)/2$ was required, I\ arbitrarily chose $a=1/12$\ and $b=7/12$. \ 

\begin{quotation}
\noindent Perform exactly the same simulation as before, except assume
treatment lengths (in minutes) are uniformly distributed on the interval
$[5,35]$. \ Less is known about this scenario than the preceding (with
exponential service times).
\end{quotation}

\noindent The thought of choosing $(a,b)=(0,2/3)$ or $(1/3,1/3)$ did not occur
to me until later. \ I\ had imagined that numerical inversion of Laplace
transforms was the only avenue available to reliably estimate the mode and median.

Just as $(a+b)/2$ is the first moment of treatment lengths,
\[%
\begin{array}
[c]{ccc}%
\xi=\dfrac{a^{2}+a\,b+b^{2}}{3}, &  & \eta=\dfrac{(a+b)\left(  a^{2}%
+b^{2}\right)  }{4}%
\end{array}
\]
are the corresponding second and third moments. \ Again, favoring pristine
theory over experiment,%
\[
F(s)=\frac{(1-\rho)s}{s-\lambda+\lambda\tfrac{\exp(-a\,s)-\exp(-b\,s)}%
{(b-a)s}};
\]
consequently \cite{Tij2-que, ZHG-que}%
\[%
\begin{array}
[c]{ccc}%
\text{mean}=-F^{\prime}(0)=\dfrac{\lambda\,\xi}{2(1-\rho)}, &  &
\text{variance}=F^{\prime\prime}(0)-F^{\prime}(0)^{2}=\dfrac{\lambda\,\eta
}{3(1-\rho)}+\dfrac{\lambda^{2}\xi^{2}}{4(1-\rho)^{2}}%
\end{array}
\]
giving $\frac{19}{48}$ and $\frac{1883}{6912}$ respectively. \ The mode
(location of the density maximum, excluding $0$) occurs when%
\[
f_{2}^{\prime}(x_{0})=0
\]
i.e.,%
\[
x_{0}=\tfrac{1}{6}\left(  1+3e^{1/6}-\sqrt{3e^{1/6}\left(  7-3e^{1/6}\right)
}\right)  =0.17405980...;
\]
solving the equation
\[%
{\displaystyle\int\limits_{c_{0}}^{c_{1}}}
f_{0}(t)dt+%
{\displaystyle\int\limits_{c_{1}}^{c_{2}}}
f_{1}(t)dt+%
{\displaystyle\int\limits_{c_{2}}^{x}}
f_{2}(t)dt=\tfrac{1}{2}%
\]
gives the median to be $0.21673428...$.

My classroom example deviates from the direction of research \cite{SBFGM-que},
which emphasizes heavy-tailed service time distributions. \ Ramsay
\cite{Ram-que}\ discovered a remarkably compact formula -- a single integral
of a non-oscillating function over the real line -- in connection with Pareto
service times. \ My formulas for uniform service times are sprawling by comparison.

The recursive solution of delay-differential equations is certainly not new
\cite{Fin0-que} but application of such to queueing theory does not seem to be
widespread. \ Counterexamples include \cite{Vla-que, BoxV-que}; surely there
are more that I've missed. \ Symbolics are mentioned in \cite{SBFGM-que}. \ 

The waiting time probability for M/D/1 seems to decay exponentially (as do the
other two cases). \ The Cram\'{e}r-Lundberg approximation \cite{Sak-que} is
applicable; alternatively, we have asymptotics \cite{Tij2-que, TSt-que}%
\[%
\begin{array}
[c]{ccc}%
\operatorname*{P}\left\{  W_{\text{que}}>x\right\}  \sim\dfrac{1-\rho}%
{\tau\,\rho-1}\exp\left[  -\lambda(\tau-1)x\right]  &  & \text{as
}x\rightarrow\infty\text{ }%
\end{array}
\]
where $\tau>1$ is the unique root of $\tau\exp\left[  -\rho(\tau-1)\right]
=1$.

Returning finally to M/U/1, let $L_{\text{sys}}$ denote the number of patients
in the system (both queue and service). \ Define%
\[%
\begin{array}
[c]{ccc}%
\tilde{F}(z)=\dfrac{(1-\rho)(1-z)\Theta(\lambda(1-z))}{\Theta(\lambda
(1-z))-z}, &  & \sigma^{2}=\xi-\dfrac{1}{\mu^{2}}.
\end{array}
\]
Under equilibrium, with $\lambda=2$, $\mu=3$ and $a=\frac{1}{12}$, $b=\frac
{7}{12}$, we have
\[
\tilde{f}(\ell)=\mathbb{P}\left\{  L_{\text{sys}}=\ell\right\}  =\frac{1}%
{\ell!}\left.  \frac{d^{\ell}\tilde{F}}{dz^{\ell}}\right\vert _{z=0}=\left\{
\begin{array}
[c]{lll}%
\frac{1}{3}=0.333333... &  & \text{if }\ell=0\\
\frac{1-e+e^{7/6}}{-3+3e}=0.289628... &  & \text{if }\ell=1\\
0.177042... &  & \text{if }\ell=2\\
0.096164... &  & \text{if }\ell=3\\
0.050209... &  & \text{if }\ell=4\\
0.025950... &  & \text{if }\ell=5
\end{array}
\right.  ;
\]
consequently \cite{CGA-que, Grif-que}%
\[
\text{mean}=\tilde{F}^{\prime}(1)=\rho+\dfrac{\rho^{2}+\lambda^{2}\sigma^{2}%
}{2(1-\rho)},
\]%
\begin{align*}
\text{variance}  &  =\tilde{F}^{\prime\prime}(1)+\tilde{F}^{\prime}%
(1)-\tilde{F}^{\prime}(1)^{2}\\
&  =\rho(1-\rho)+\dfrac{(3-2\rho)(\rho^{2}+\lambda^{2}\sigma^{2})}{2(1-\rho
)}+\dfrac{(\rho^{2}+\lambda^{2}\sigma^{2})^{2}}{4(1-\rho)}+\frac{\lambda
^{3}\eta}{3(1-\rho)}+\frac{\lambda^{4}\rho\,\xi^{2}}{4(1-\rho)^{2}}%
\end{align*}
giving $\frac{35}{24}$ and $\frac{4547}{1728}$ respectively. \ The final term
for the variance is missing in \cite{Grif-que}. \ The probability generating
function $\tilde{f}(\ell)$ seems to decay geometrically, but details
surrounding the exact limit of successive ratios have not been verified
\cite{CGA-que}.

I am grateful to innumerable software developers, as my \textquotedblleft
effective\textquotedblright\ formulas are too lengthy to be studied in any
traditional sense. \ Mathematica routines NDSolve for DDEs and
InverseLaplaceTransform (for Mma version $\geq12.2$) plus ILTCME
\cite{CME-que} assisted in numerically confirming many results. \ R
steadfastly remains my favorite statistical programming language. \ A student
asked in 2006 for my help in writing a relevant R\ simulation, leading to the
computational exercises described here and to my abiding interest in queues
\cite{Fi1-que, Fi2-que}.%
\begin{figure}
[ptb]
\begin{center}
\includegraphics[
height=4.1511in,
width=6.4083in
]%
{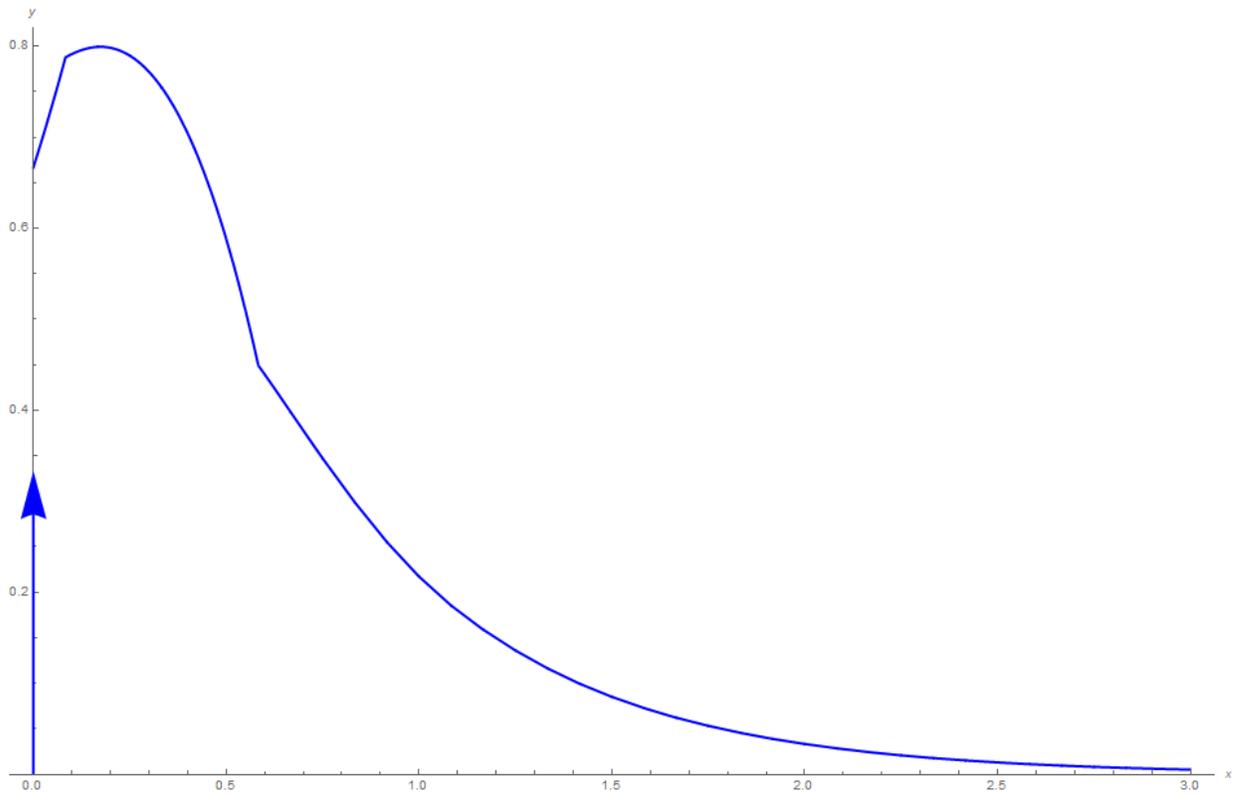}%
\caption{Waiting time density plot $y=f(x)$ for Uniform[$\frac{1}{12},\frac
{7}{12}$] service.}%
\end{center}
\end{figure}
\begin{figure}
[ptb]
\begin{center}
\includegraphics[
height=4.2186in,
width=6.5484in
]%
{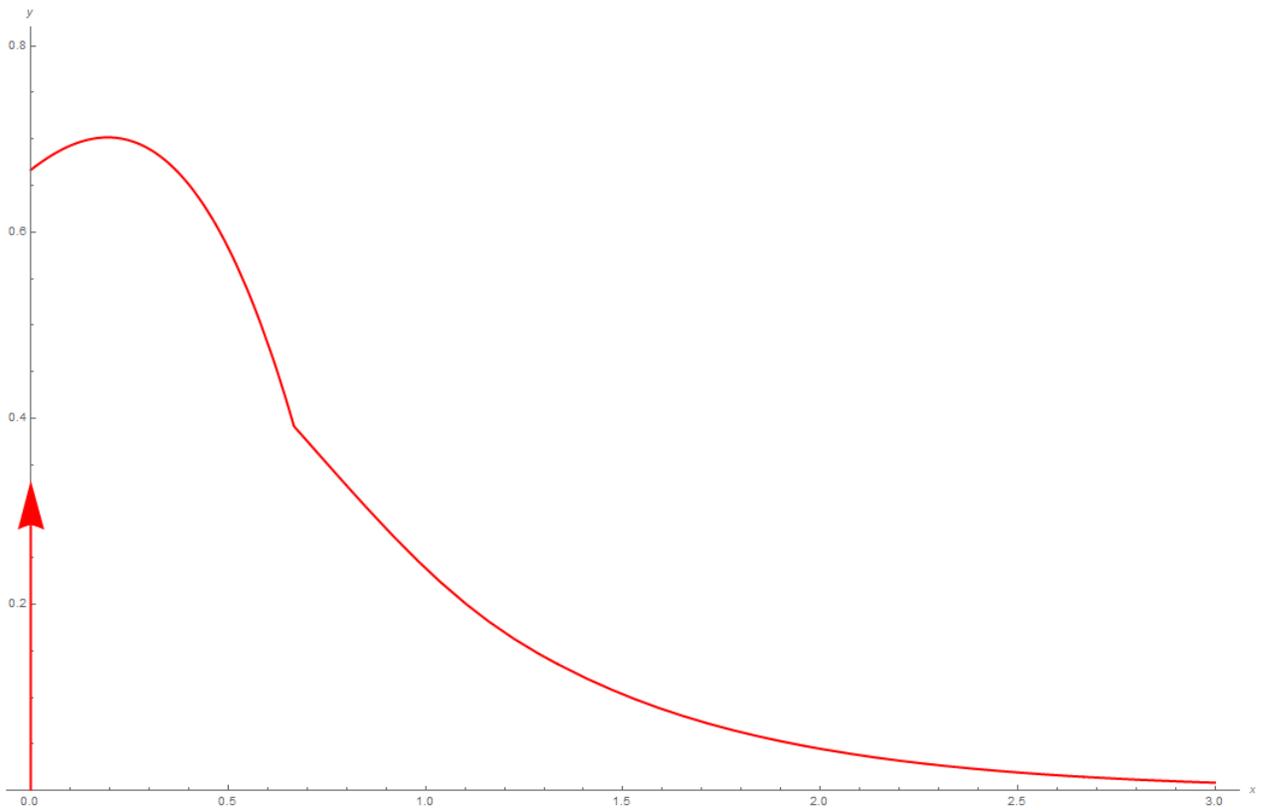}%
\caption{Waiting time density plot $y=g(x)$ for Uniform[$0,\frac{2}{3}$]
service.}%
\end{center}
\end{figure}
\begin{figure}
[ptb]
\begin{center}
\includegraphics[
height=4.3145in,
width=6.6928in
]%
{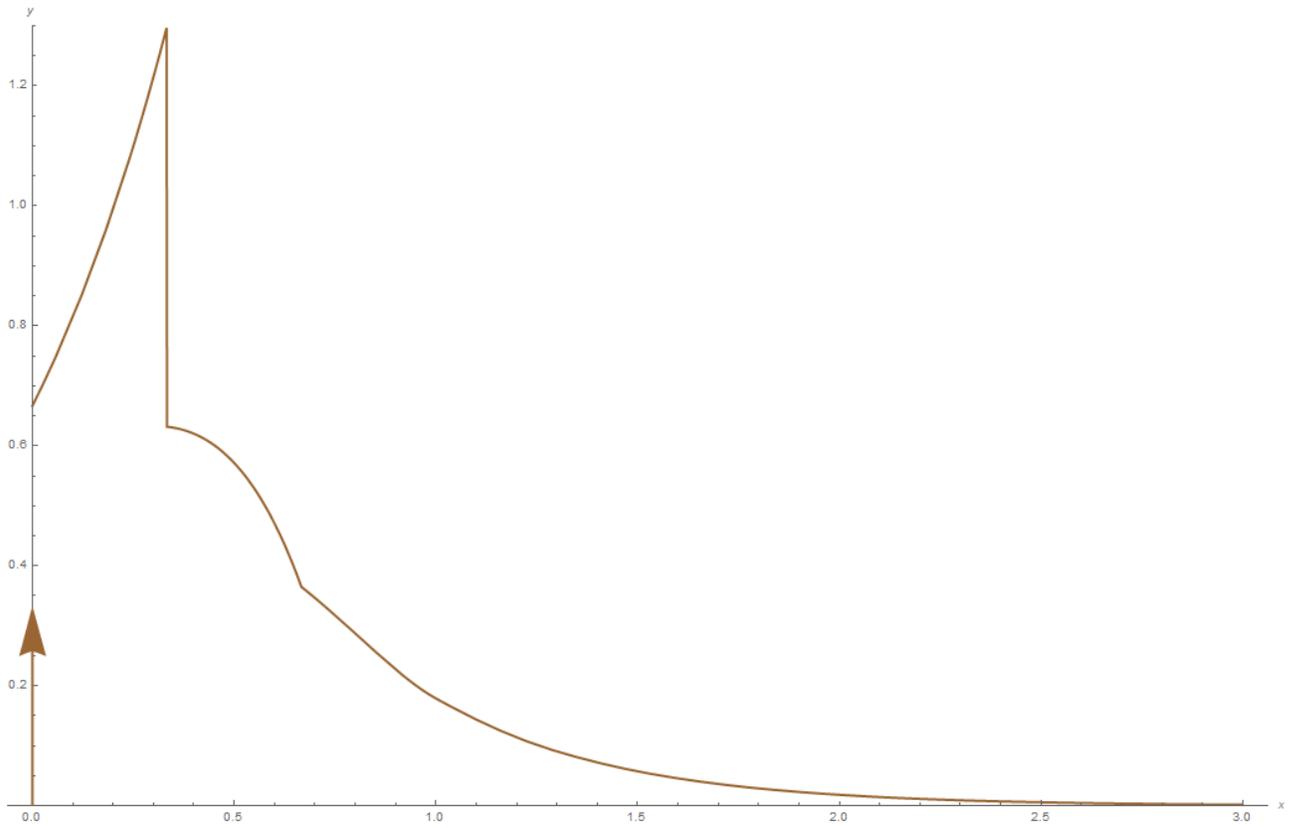}%
\caption{Waiting time density plot $y=h(x)$ for Deterministic[$\frac{1}{3}$]
service.}%
\end{center}
\end{figure}

\section{Addendum}

For completeness, another M/D/1 result is provided. Under equilibrium, with
$\lambda=2$, $\mu=3$ and $a=\frac{1}{3}=b$, we have
\[
\tilde{h}(\ell)=\mathbb{P}\left\{  L_{\text{sys}}=\ell\right\}  =\left\{
\begin{array}
[c]{lll}%
\frac{1}{3}=0.333333... &  & \text{if }\ell=0\\
\frac{-1+e^{2/3}}{3}=0.315911... &  & \text{if }\ell=1\\
\frac{-5e^{2/3}+3e^{4/3}}{9}=0.182481... &  & \text{if }\ell=2\\
\frac{8e^{2/3}-21e^{4/3}+9e^{2}}{27}=0.089494... &  & \text{if }\ell=3\\
\frac{-22e^{2/3}+180e^{4/3}-243e^{2}+81e^{8/3}}{243}=0.042035... &  & \text{if
}\ell=4\\
\frac{14e^{2/3}-312e^{4/3}+972e^{2}-891e^{8/3}+243e^{10/3}}{729}=0.019607... &
& \text{if }\ell=5
\end{array}
\right.  ;
\]
giving $\frac{4}{3}$ and $\frac{56}{27}$ for the mean and variance, yet
unverified limit of successive ratios.

\end{document}